\documentclass[911pt]{article}
\usepackage{amsmath,amsfonts,latexsym,amssymb,amsthm}
\usepackage{url}

\newcommand{\Q}{\mathbb{Q}}

\newcommand{\Z}{\mathbb{Z}}

\newcommand{\F}{\mathbb{F}}
\newcommand{\PP}{\mathbb{P}}

\begin{document}

\title{Elliptic curves with large torsion and positive rank over number fields of small degree
and ECM factorization}
\author{Andrej Dujella and Filip Najman}
\date{}
\maketitle
\begin{abstract}
In this paper, we present several methods for construction of elliptic curves with large torsion group and positive rank over number fields of small degree. We also discuss potential applications of such
curves in the elliptic curve factorization method (ECM).
\end{abstract}

\section{Introduction.}

Let $E$ be an elliptic curve over $ \mathbb{Q}$.
By the Mordell-Weil theorem, the group $E( \mathbb{Q})$ of rational points
on $E$ is a finitely generated abelian group. Hence, it is the
product of the torsion group and $r\geq 0$ copies of an infinite cyclic
group: $ E( \mathbb{Q}) \cong E( \mathbb{Q})_{\rm
tors} \times { \mathbb{Z} }^r$.
By Mazur's theorem, we know that $E( \mathbb{Q})_{\rm tors}$ is one
of the following 15 groups:
$ \mathbb{Z}/n\mathbb{Z}$ with $1\leq n \leq 10$ or $n=12$,
$ \mathbb{Z}/2\mathbb{Z} \times \mathbb{Z}/2m\mathbb{Z}$ with $1\leq
m\leq 4$.

It is well-known that elliptic curves have applications in cryptography and also in factorization of large integers and primality proving. The main idea is to replace the group $\mathbb{F}_p^*$ with (fixed) order $p-1$, by the group $E(\mathbb{F}_p)$ with more flexible order. Namely, by the Hasse theorem we have
$$ p+1-2\sqrt{p} < |E(\mathbb{F}_p)| < p+1+2\sqrt{p}. $$

In particular, elliptic curves with large torsion and positive rank (it is necessary for an elliptic curve to have positive rank for it to be used for factorization) over the rationals have long been used for factorization, starting with Montgomery, Atkin and Morain (\cite{Mont,A-M}).
We say that an integer $m$ is \emph{$n$-smooth}, for some fixed value $n$ if
all the prime divisors of $m$ are less or equal than $n$.
Choosing elliptic curves $E$ for the elliptic curve factoring method (ECM),
one wants to choose elliptic curves such that the order $E(\mathbb F_p)$ is smooth (for more details about ECM, see \cite{len}, where the method was introduced).

Standard heuristics say that larger torsion of $E(\mathbb Q)$ implies a greater probability that $|E(\mathbb F_p)|$ is smooth. This is because the torsion of $E(\mathbb Q)$ will inject into $E(\mathbb F_p)$ for all primes $p$ of good reduction, making $|E(\mathbb F_p)|$ divisible by the order of the torsion of $E(\mathbb Q)$. But this is not necessary so straightforward, as shown in \cite{bbdn}, as a curve with smaller $E(\mathbb Q)_{{\rm tors}}$ can have much larger torsion over fields of small degree, giving all together a greater probability of $|E(\mathbb F_p)|$ to be smooth.

As we shown in Section \ref{sec:exp} of this paper, this is especially true in some instances, when for some reasons the properties of the prime factors of the numbers that are to be factored are known. One can then do better by choosing elliptic curves with large torsion and positive rank over some small extension of $\Q$. Nice explicit examples of this approach, for factoring large numbers (Cunningham numbers in this case) by using elliptic curves over number fields, have been provided recently by Brier and Clavier \cite{bc}. They used elliptic curves over cyclotomic fields with torsion groups $\Z/3\Z \oplus \Z/6\Z$ and $\Z/4\Z \oplus \Z/4\Z$.

Also, they tried to construct elliptic curves over cyclotomic fields with torsion $\Z/5\Z \oplus \Z/5\Z$ and $\Z/4\Z \oplus \Z/8\Z$ and positive rank, but failed. Recently, examples of such curves have been found in \cite{bbdn}. In this paper we show how to construct such curves systematically in Section \ref{sec:kons}.

It is very useful to have families of curves suitable for use in ECM, as when one curve fails to yield a factorization, another one can be used. For the largest torsions from Mazur's theorem, families of elliptic curves with large torsion and positive rank can be found in \cite{A-M}. Each of these is a parameterized family of elliptic curves, in other words there exists a surjective map from $\PP^1(\Q)$ (minus some points) to this family of elliptic curves. This can also be understood as a map from a genus 0 curve to the family of elliptic curves.

This is the best possible case, but for some torsion groups, such maps are yet unknown, and in most cases theoretically impossible. In these cases we can instead construct a map from some curve of positive genus to the family of elliptic curves with large torsion and positive rank. An example of such a map was constructed in \cite{bc}, where elliptic curves with torsion $\Z /3\Z \oplus \Z /6\Z$ and positive rank were parameterized by an elliptic curve of positive rank. Note that it is preferable for the genus of the curve to be as small as possible. In particular if we have a map from an elliptic curve with positive rank, this again allows us to construct infinitely many curves in such a manner. In Section \ref{sec:duje} we construct 2 new examples of families of elliptic curves with positive rank and torsion isomorphic to $\Z /3\Z \oplus \Z /6\Z$ over $\Q(\sqrt{-3})$, that are each parameterized by an elliptic curve of positive rank. We also construct a family of elliptic curves with torsion $\Z /5\Z \oplus \Z /5\Z$ and positive rank over $\Q(\zeta_5)$ that is parameterized by a genus 2 curve.

Families of elliptic curves with high rank and prescribed torsion over the rationals have been of great interest historically (see \cite{tors2} for a list of references). In recent years there has been also a great interest for such families over quadratic fields \cite{Rab, mir1, dm}. In Section \ref{sec:duje} we also construct a parameterized family of elliptic curves with torsion $\Z /3\Z \oplus \Z /3\Z$ and rank 2, which was previously unknown. The benefit of using rank 2 curves instead of rank 1 curves for ECM might be in the case when the reduction of one generator is not of smooth order. With a rank 2 curve one can hope that the other generator might reduce to a point of smooth order, although the probability for this is not very high, and thus the algorithmical benefits of rank $2$ are not clear at the moment.

When we have a map from a curve of genus $>1$, then we can in this way always construct only finitely many curves with the desired property. In these cases it is actually more useful to have a procedure to construct such curves. In Section \ref{sec:kons} we give a procedure that constructs elliptic curves with torsion groups $\Z /4\Z \oplus \Z /8\Z$, $\Z /5\Z \oplus \Z /5\Z$ and $\Z /6\Z \oplus \Z /6\Z$ and positive rank over quartic fields. We also give some examples of curves constructed in this way.

\section{Choosing curves for ECM depending on the primes}
\label{sec:exp}
In this section, we show how the chance for an elliptic curve $E$ over $\Q$ to be smooth over $\F_p$ really depends on the set of primes that $p$ runs through. The ideas presented here expand on the arguments outlined in \cite{bbdn}. We will show that knowing the splitting behavior of the primes in some extension of $\Q$ over which $E$ has large torsion can be used to determine how likely $E$ is to be smooth.

Let $E$ be an elliptic curve with torsion $T_1$ over $\Q$ and $T_2$ over some number field $K$ of relatively small degree and suppose $|T_1|<|T_2|$. Suppose for simplicity that $K$ is Galois (one could make the same argument with non-Galois extensions, but it would be a bit more messy). Then one can see if the rational prime $p$ splits completely in $K$, then the whole torsion of $E(K)$ will inject into $E(\F_p)$ (see  \cite[Proposition 3.1, pp.~176]{sil}). Thus we expect $E(\F_p)$ to be smooth more often when $p$ runs through the primes that split completely in $K$ than over the set of all primes.

We test this heuristic by choosing the following 8 elliptic curves:
$$E_0: y^2 = x^3 + 3,$$
$$E_7: y^2 - 55xy - 448y = x^3 - 448x^2,$$
$$E_9: y^2 - 47xy - 624y = x^3 - 624x^2,$$
$$E_{12}:  y^2 + 19/40xy - 273/400y = x^3 - 273/400x^2,$$
$$E_{2\times 8}: y^2 = x^3 + 54721/225x^2 + 4096x,$$
$$E_{4\times 8}: y^2 + xy + y = x^3 + x^2 - 52431x - 2731947,$$
$$E_{5 \times 5}:  y^2 + y = x^3 + x^2 - 5092900x + 709824595630,$$
$$E_{6 \times 6}:  y^2 + xy + y = x^3 - 371066x - 47384980.$$

The curves were chosen in the following way: $E_0$, $E_7$, $E_9$ and $E_{12}$ have trivial torsion, 7-torsion, 9-torsion and 12-torsion over $\Q$, respectively. The curve $E_{2\times 8}$ has torsion $\Z /2\Z \oplus \Z /8\Z$ over $\Q$. The curves $E_{4\times 8}, E_{5 \times 5}$ and $E_{6 \times 6}$ were chosen such that they have torsions $\Z /4\Z \oplus \Z /8\Z$, $\Z /5\Z \oplus \Z /5\Z$ and $\Z /6\Z \oplus \Z /6\Z$ over some biquadratic field. Note that $E_{4\times 8}$ has torsion $\Z /2\Z \oplus \Z /2\Z$ over $\Q$, $E_{5\times 5}$ has 5-torison over $\Q$ and $E_{6 \times 6}$ has 6-torsion over $\Q$. All the listed elliptic curves have rank 1 over $\Q$, so the rank should not influence the smoothness results.
The curve $E_{5 \times 5}$ is taken from \cite{bbdn}. Note that new curves with the same property
will be constructed in Section \ref{sec:duje} and \ref{sec:kons}.

We now test how often the reductions of these curves $|E(\F_p))|$ have 100-smooth order, depending on the set of primes that we choose. We will test subsets of the set $\{p_n| 50 \leq n \leq 10050\}$, where $p_n$ denotes the $n$-th prime (we choose $n\geq 50$ to get rid of the primes of bad reduction).
Let
$$ A=\{p_n| 50 \leq n \leq 10050\},$$
$$ B=\{p\in A | \left( \frac{-143}{p }\right)=1\},$$
$$ C=\{p\in A | \left(\frac{-143}{p}\right)=-1\},$$
$$ D=\{p\in A | p\equiv 1\pmod 5\},$$
$$ E=\{p\in A | \left(\frac{-3}{p}\right)=1 \text{ and } \left(\frac{217}{p}\right)=1\},$$
$$ F=\{p\in A | \left(\frac{-1}{p}\right)=1 \text{ and } \left(\frac{-7}{p}\right)=1\}.$$

The elliptic curve $E_{12}$ has torsion $\Z /2\Z \oplus \Z /12\Z$ over $\Q(\sqrt{-143})$, and $ \Z /12\Z$ over all other quadratic fields, so we expect this curve to be more likely smooth over $\F_p$, where $p$ splits in $\Q(\sqrt{-143})$, then when $p$ does not split. We compare this by examining the reductions over the primes from the sets $B$ and $C$. The curve $E_{5\times 5}$ has torsion $\Z /5\Z \oplus \Z /5\Z$ over $\Q(\zeta_5)$. Recall that a rational prime $p$ splits completely in $\Q(\zeta_5)$ if and only if $p\equiv 1 \pmod 5$. Thus we expect $|E_{5\times 5}(\F_p)|$ to be more likely smooth when $p\equiv 1 \pmod 5$. The curves $E_{4\times 8}$ and $E_{6\times 6}$ have torsions $\Z /4\Z \oplus \Z /8\Z$, and $\Z /6\Z \oplus \Z /6\Z$ over the fields $\Q(\sqrt{-1},\sqrt{-7})$ and $\Q(\sqrt{-3},\sqrt{217})$, so the sets $E$ and $F$ have been chosen in a way such that we expect $|E_{4\times 8}(\F_p)|$ and $|E_{6\times 6}(\F_p)|$, respectively, to have a greater probability to be smooth.

In the table below we list for each set $S$ and each curve $E$ the number of 100-smooth values of $|E(\F_p)|$, where $p\in S$.

\begin{center}
\begin{tabular}{|c|c|c|c|c|c|c|}
\hline
el. curve & $A$ & $B$ & $C$ & $D$ & $E$ & $F$ \\
\hline
$E_0$ & 2822 & 1453 & 1369 & 643 & 522 & 633 \\
\hline
$E_7$ & 4275 & 2115 & 2160 & 1020 & 1014 & 1066 \\
\hline
$E_9$ & 4635 & 2306 & 2329 & 1110 & 1226 & 1125\\
\hline
$E_{12}$      & 5133 & 2852 & 2281 & 1290 & 1302 & 1288\\
\hline
$E_{2\times 8}$      & 5110 & 2587 & 2523 & 1245 & 1206 & 1295\\
\hline
$E_{4\times 8}$      & 4317 & 2141 & 2176 & 1059 & 1098 & 1440\\
\hline
$E_{5\times 5} $     & 4376 & 2137 & 2239 & 1448 & 1047 & 1074\\
\hline
$E_{6\times 6}  $    & 4817 & 2396 & 2421 & 1201 & 1505 & 1138\\
\hline

\end{tabular}
\end{center}

We see that when we run through all primes (the set $A$), then the curves most likely to be smooth are $E_{12}$ and $E_{2 \times 8}$, with them approximately  being equally likely to be smooth. In the columns with the sets $B$ and $C$ we see that $E_{12}$ is approximately $10\%$ more likely to be smooth than $E_{2 \times 8}$  when reducing modulo the primes from $B$, and more than $10\%$ less likely to be smooth when reducing modulo the primes from $C$. We also see that the curves $E_{4\times 8}, E_{5 \times 5}$ and $E_{6 \times 6}$ are most likely to be smooth when reducing modulo the primes from $D$, $E$ and $F$, respectively, outperforming all the other curves quite convincingly in each case.

These result strongly suggest that, when performing ECM, if one knows the splitting behavior of the primes over which the curves are going to be reduced, one can in some cases do considerably better than just choosing elliptic curves with the largest possible torsion over $\Q$.

\section{A method for finding subfamilies with larger rank.}
\label{sec:duje}
In this section we will describe a method for construction of families of elliptic curves
with certain property and relative large rank. We assume that a family of elliptic curves
with that property is known, and we show how to find its subfamily with larger generic rank.

Let $\mathbb{K}$ be a number field, $E \,:\, y^2=f(x)$ an elliptic curve
over $\mathbb{K}(T)$ and let $\Delta$ be its
discriminant. Assume that $E$ has a nontrivial torsion group. We seek congruences
of the shape $x\equiv x_0 \pmod{\delta}$, where $\delta$ is a factor of $\Delta$,
which are satisfied by $x$-coordinates of some of the torsion points
(and other known points on $E(\mathbb{K}(T))$ in the case of curves with
positive generic rank and known points of infinite order).
Then we search for further (nontorsion) points on $E$ with $x$-coordinate of the form
$x=x_0+\delta p$, where $p$ is a polynomial over $\mathbb{K}$ with small degree (say
$\deg p \leq 2$). We insert  $x=x_0+\delta (aT^2+bT+c)$ into $f(x)$, get rid of the quadratic factor,
and impose the condition that the discriminant of the remaining polynomial in $T$ is $0$.
This gives us several equations for $a,b,c$. Substituting the obtained conditions, we repeat
the procedure (getting rid of quadratic factors and asking that the remaining polynomial in $T$
has zero discriminant). Finally, we get a condition of the form $g(T)=z^2$. If the condition
corresponds to a curve of genus $0$ or $1$ with a $\mathbb{K}$-rational point, then we obtain
a subfamily for $E$ with potentially larger rank. That the rank indeed increases, can be checked
by finding a suitable specialization for which the corresponding points are independent.

A variant of this procedure has been previously successfully used for finding generators
of some high rank elliptic curves over $\mathbb{Q}$ with relatively large torsion group
which can be found in \cite{tors}. A motivation for this methods comes from the Lutz-Nagell
theorem which says that the torsion points, but possibly also some other integer points,
satisfy $y\equiv 0\pmod{\delta}$ for a factor $\delta$ of $\Delta$. But this implies
$x\equiv x_0 \pmod{\delta}$ for $x_0$ from a finite set.

We will illustrate the above method by constructing the first known example of
an elliptic over $\mathbb{Q}(\sqrt{-3})(T)$ with torsion group $\mathbb{Z}/3\mathbb{Z} \oplus \mathbb{Z}/3\mathbb{Z}$ and rank $\geq 2$. We will also apply
variants of this method to curves with torsion groups $\mathbb{Z}/3\mathbb{Z} \oplus
 \mathbb{Z}/6\mathbb{Z}$ and $\mathbb{Z}/5\mathbb{Z} \oplus \mathbb{Z}/5\mathbb{Z}$,
and we will obtain results comparable to those from \cite{bc} and \cite{bbdn}.

\bigskip

\framebox{$\mathbb{Z}/3\mathbb{Z} \oplus \mathbb{Z}/3\mathbb{Z}$}
\smallskip

Our starting point is the curve over $\mathbb{Q}(\sqrt{-3})(T)$ with
torsion group $\mathbb{Z}/3\mathbb{Z} \oplus \mathbb{Z}/3\mathbb{Z}$ and rank
$\geq 1$ found by Rabarison \cite{Rab}:
$$ Y^2 = X^3+(108+T^6)X^2+(144T^6+3888)X+64T^12+3456T^6+46656. $$
Its discriminant is
$\Delta=-4096T^{12}(T^2+3)^3(T^2+3T+3)^3(T^2-3T+3)^3$.
The torsion points over $\mathbb{Q}(T)$ are $\mathcal{O}$, $(-4T^4+12T^2-36, -4T^7+12T^5-36T^3)$,
$(-4T^4+12T^2-36, 4T^7-12T^5+36T^3)$ (with an additional
point $(4/3T^6-36, 4\sqrt{-3}/9T^9+12\sqrt{-3}T^3)$ of order $3$ over $\mathbb{Q}(\sqrt{-3})(T)$).
It has positive rank, with the point $(0, 8T^6+216)$ of infinite order.

Applying the above procedure to the points with the first coordinate of the form
$X=(aT+b)(T^2+3T+3)$, leads to
$X=(12T-12)(T^2+3T+3)$ and the condition
$13T^2-42T+93 = z^2$, which is a genus $0$ curve with the parametrization
$$ T=\frac{u^2-16u-29}{u^2-13}. $$
By taking $u=1$, we get the curve
$$ Y^2=X^3+1850293/729X^2+28659904/81X+205347524322304/531441 $$
and the independent points of infinite order
$(0,14329952/729)$ and $(7904/9$, $42080896/729)$.
Hence, we have constructed a curve over $\mathbb{Q}(u)$ with rank $\geq 2$ and with
torsion group over $\mathbb{Q}(\sqrt{-3})(u)$ isomorphic to $\mathbb{Z}/3\mathbb{Z} \oplus \mathbb{Z}/3\mathbb{Z}$.
Explicitly, the curve is
\vspace{-2ex}

{\footnotesize
\begin{eqnarray*}
y^2 &=& x^3 + (1116118693+2352294282u^2+508999635u^4-26095764u^6+823995u^8\\
& &  \mbox{}-4758u^{10}+109u^{12}+1969070304u+
1658452000u^3-818496u^5+28224u^7-68000u^9 \\
& & \mbox{}-96u^{11}) x^2+
576(7u^2+32u+37)(u^2-16u+67)(u^4-8u^3+30u^2+232u+337)\\
& & \mbox{}\times (u^4+16u^3+246u^2+304u+217)(u^2-13)^6 x \\
& & \mbox{}+
1024(7u^2+32u+37)^2(u^2-16u+67)^2(u^4-8u^3+30u^2+232u+337)^2\\
& & \mbox{}\times(u^4+16u^3+246u^2+304u+217)^2(u^2-13)^{6},
\end{eqnarray*}}%
with independent points of infinite order
\begin{eqnarray*}
P_1 &=& (0,
32(7u^2+32u+37)(u^2-16u+67)(u^4-8u^3+30u^2+232u+337) \\
& & \mbox{}\times (u^4+16u^3+246u^2+304u+217)(u^2-13)^3),
\end{eqnarray*}
\begin{eqnarray*}
P_2 &=& (-192(u+1)(7u^2+32u+37)(u^2-16u+67)(u^2-13)^3,32\\
& &
(7u^2+32u+37) (u^2-16u+67)(u^2+2u+13)(u^2-16u-29)^3(u^2-13)^3).
\end{eqnarray*}

\newpage

\framebox{$\mathbb{Z}/3\mathbb{Z} \oplus \mathbb{Z}/6\mathbb{Z}$}
\smallskip

Since in this case we have a point of order $2$, we can write the curve in the form
$y^2 = x^3+ax^2+bx$, and we get the general form of curves with torsion group
$\mathbb{Z}/3\mathbb{Z} \oplus \mathbb{Z}/6\mathbb{Z}$ over $\mathbb{Q}(\sqrt{-3})$
by taking
$A=-24(3T^2-4)(9T^4-36T^3+72T^2-48T+16)$,
$B=144(T-2)^3(3T-2)^3(3T^2+4)^3$ (see \cite{Rab,mir1}).
The discriminant is
$\Delta = 97844723712T^3(3T^2-6T+4)^3(T-2)^6(3T-2)^6(3T^2+4)^6$.
We search for points on the curve such that their $x$-coordinate is
a factor $b_1$ of $B$. (A similar method was used in \cite{dm}
for finding a curve over $\mathbb{Q}(i)(T)$ with torsion group $\mathbb{Z}/4\mathbb{Z} \oplus \mathbb{Z}/4\mathbb{Z}$ and rank $\geq 2$). This leads to the condition that
$b_1 + A + B/b_1$ is a perfect square.
Several such factors lead to the condition which correspond
to curve of genus 1.
E.g.
\begin{itemize}

\item[(i)] for
$b_1=-12(T-2)^3(3T^2+4)(3T-2)$ we get
the genus $1$ curve $-27T^4+72T^3-144T^2+96T=z^2$ with rank $1$ (the minimal Weierstrass equation is
$Y^2+Y = X^3-34$);

\item[(ii)] for
$b_1=3/2(T-2)^2(3T^2+4)^2(3T-2)^3$ we get
the genus $1$ curve
$162T^3-324T^2+216T-1200=z^2$ with rank $2$ (the minimal Weierstrass equation is
$Y^2 = X^3-648$);

\item[(iii)] for
$b_1=4(T-2)(3T^2+4)(3T-2)^3$ we get
$-6T^3+12T^2-8T+4=z^2$ with rank $1$ (the minimal Weierstrass equation is
$Y^2+Y = X^3+1$).
\end{itemize}
The example (i) is equivalent to the example given in \cite{bc}, while the examples (ii) and (iii)
give new examples of infinite families of curves over $\mathbb{Q}(\sqrt{-3})$ with
torsion group $\mathbb{Z}/3\mathbb{Z} \oplus \mathbb{Z}/6\mathbb{Z}$ and positive rank.

\bigskip

\framebox{$\mathbb{Z}/5\mathbb{Z} \oplus \mathbb{Z}/5\mathbb{Z}$}
\smallskip

The curve $y^2 = x^3 + Ax +B$, where
$A=-27(T^{20}+228T^{15}+494T^{10}-228T^5+1)$,
$B=54(T^{30}-522T^{25}-10005T^{20}-10005T^{10}+522T^5+1)$
has torsion group $\mathbb{Z}/5\mathbb{Z} \oplus \mathbb{Z}/5\mathbb{Z}$ over cyclotomic field
$\mathbb{Q}(\zeta_5)$ (see \cite{bc}). Its discriminant is
$\Delta=2176782336T^5(T^2-T-1)^5(T^4-2T^3+4T^2-3T+1)^5(T^4+3T^3+4T^2+2T+1)^5$.
The $x$-coordinates of the torsion points over $\mathbb{Q}(T)$,
$3T^{10}+36T^9+72T^8+108T^7+180T^6+90T^5+72T^4-36T^3+36T^2+3$
and
$3T^{10}+36T^8+36T^7+72T^6-90T^5+180T^4-108T^3+72T^2-36T+3$,
satisfy the congruence
\begin{eqnarray*}
x &\equiv& 33T(T^2-T+1)(T^5+2T^4+3T^3-T^2+T-1)+3 \\
 & & \pmod{T(T^4-2T^3+4T^2-3T+1)(T^4+3T^3+4T^2+2T+1)}.
\end{eqnarray*}
Searching for nontorsion points of the form
$x=(aT^2+bT+c)(T^4-2T^3+4T^2-3T+1)(T^4+3T^3+4T^2+2T+1)+ 33T(T^2-T+1)(T^5+2T^4+3T^3-T^2+T-1)+3$,
we were not able to reach the condition leading to curves of genus $0$ or $1$.
However, by taking $aT^2+bT+c = -6T^2+6T-27$,
we obtain the condition $T^5-18 = z^2$ which gives a curve of genus $2$.
By taking $T=3$, we get the curve
$$ y^2+y = x^3+x^2-226248x-20170186 $$
with rank $1$ and torsion $\mathbb{Z}/5\mathbb{Z} \oplus \mathbb{Z}/5\mathbb{Z}$.
The point $(-132, 2722)$ is of infinite order.

\section{Constructing individual curves with torsion \\ $\Z /4\Z \oplus \Z /8\Z$, $\Z /5\Z \oplus \Z /5\Z$ and $\Z /6\Z \oplus \Z /6\Z$}
\label{sec:kons}

As mentioned in the introduction, we often cannot construct a surjective map from a curve of genus $\leq 1$ to a family of elliptic curves over $\Q$ with prescribed torsion and positive rank. This is exactly the case for the torsions $\Z /4\Z \oplus \Z /8\Z$, $\Z /5\Z \oplus \Z /5\Z$ and $\Z /6\Z \oplus \Z /6\Z$, which are interesting for ECM, as shown in Section \ref{sec:exp}.

Curves with torsion $\Z /6\Z \oplus \Z /6\Z$ over a biquadratic field are the easiest case and we deal with them first. In fact one can use the elliptic curves over $\Q(\sqrt{-3})$ obtained in Section \ref{sec:duje}, and then for each curve construct an elliptic curve with torsion $\Z /6\Z \oplus \Z /6\Z$ over the quadratic extension of $\Q(\sqrt{-3})$ obtained by adjoining the root of the discriminant of the elliptic curve. Each curve constructed in this way have rank at least 2.

We are left with the torsion groups $\Z /4\Z \oplus \Z /8\Z$ and $\Z /5\Z \oplus \Z /5\Z$. We give a method of constructing curves with these torsion groups and positive rank. In theory, this can give us infinitely many elliptic curves with the desired torsion and positive rank. This method is also very useful for practical purposes, easily generating many curves with the desired properties.

The starting point is \cite{jk4}, where a method of constructing infinitely many elliptic curves with the before-mentioned torsion groups is given. We then sieve through the constructed curves in search of elliptic curves with positive rank. We do not test directly whether the curve has positive rank but instead use the fact that in this construction, for all three torsion groups, all the obtained elliptic curves are rational, and that the torsion group we desire is defined over a biquadratic field.

We then use the fact that if $K$ is a number field, $L$ a quadratic extension of $K$, $L=K(\sqrt d)$, and $E$ an elliptic curve defined over $K$, then
$$ {\rm rank}(E(L))={\rm rank}(E(K))+{\rm rank}(E^{(d)}(K)).$$
As we are interested in the rank of a rational elliptic curve $E$ over a biquadratic field $K$, one can see that the rank of $E(K)$ is the sum of the rank of the $E(\Q)$ and 3 of its twists.

To construct elliptic curves with torsion $\Z /4\Z \oplus \Z /8\Z$, we do the following:
\begin{enumerate}

\item Construct a rational elliptic curve with torsion $\Z /4\Z \oplus \Z /8\Z$ over a biquadratic field using the methods from \cite{jk4}. We start by taking $t\in \Q^*$. Define $v=(t^4-6t^2+1)/(4(t^2+1)^2$. Let $L$ be obtained by adjoining the root of $t^4-6t^2+1$ to $\Q(i)$ and let $a=v^2-1/16$. Let $E$ be the elliptic curve defined by
$$E: y^2+xy-ay=x^3-ax^2.$$
\item Check whether any of the twists of $E$ by $-1$, $t^4-6t^2+1$ and $-(t^4-6t^2+1)$ have root number -1. If they do, then the Birch-Swinnerton--Dyer conjecture suggest that the rank is odd and hence positive.

\item If any of the twists have root number -1, then search for points of infinite order on it. If none of the twists have root number -1, choose another $t$ and start over.

\item If a point of infinite order is found then we have found an elliptic curve with the desired properties! Note that it does not matter which twist we choose, as they are all isomorphic over $L$. If a point of infinite order is not found on any twist with odd root number, choose another $t$ and start over.
\end{enumerate}

 Using this procedure we can easily construct many elliptic curves with torsion $\Z /4\Z \oplus \Z /8\Z$ and positive rank over a biquadratic field. For example, the value $t=4$ gives us the elliptic curve
 $$y^2 + xy + 3600/83521y = x^3 + 3600/83521x^2$$
 with a point of infinite order $(-30/289, 3900/83521)$ and torsion $\Z /4\Z \oplus \Z /8\Z$ over $\Q(i,\sqrt{161}).$

The value $t=3$ gives the curve
$$y^2 = x^3 - 67950603/390625x – 126442451898/244140625$$
which has a point of infinite order $(-3549/625, 10584/625)$ and torsion $\Z /4\Z \oplus \Z /8\Z$ over $\Q(i,\sqrt{-7})$.

\medskip

In a only slightly different way one can construct elliptic curves with torsion $\Z /5\Z \oplus \Z /5\Z$ over $\Q(\zeta_5)$. We will not write the formulas explicitly as above, as one can find them in \cite{jk4}.

 Here the only difference is that as $\Q(\zeta_5)$ is not a biquadratic, but an extension with Galois group $\Z/4\Z$, one cannot look at 3 rational twists of the constructed elliptic curve $E$, but only at the twist by 5 (as $\Q(\sqrt 5)$ is a subfield of $\Q(\zeta_5)$). Thus we look at only 2, instead of 4 curves for root number -1.

 Nevertheless, we are able to easily construct the curves
$$y^2 = x^3 - 10605390625/10460353203x -
4238740478515625/22236242266222092$$
with  a point of infinite order $(-9875/177147,75625/3188646)$ and torsion $\Z/5\Z \oplus \Z/5\Z $ over $\Q(\zeta_5)$ and
$$y^2 = x^3 + 147734375/50331648x +
1010986328125/927712935936$$
with the point of infinite order $(15625/12288, 171875/65536)$ and torsion $\Z/5\Z \oplus \Z/5\Z $ over $\Q(\zeta_5)$. These curves were obtained by inserting the values $t=1/3$ and $t=-1/2$, respectively in the formulas in \cite{jk4}.

\bigskip

\noindent {\small \sc A. Dujella,  Department of Mathematics, University of Zagreb, \\
Bijeni\v{c}ka cesta 30, 10000 Zagreb, Croatia \\
{\tt duje@math.hr}}

\bigskip

\noindent {\small \sc F. Najman, Department of Mathematics, University of Zagreb, \\
Bijeni\v{c}ka cesta 30, 10000 Zagreb, Croatia \\
{\tt fnajman@math.hr}}


\begin{thebibliography}{x}

\bibitem{A-M}
A. O. L. Atkin and F. Morain, \emph{Finding suitable curves for the elliptic curve method of factorization}, Math. Comp. \textbf{60} (1993), 399--405.

\bibitem{bbdn}
J. Bosman, P. Bruin, A. Dujella and F. Najman, \emph{Ranks of elliptic curves with prescribed torsion over number fields}, preprint.

\bibitem{bc}
\'E. Brier, C. Clavier,
\emph{New families of ECM curves for Cunningham Numbers}. In: Proceedings of ANTS IX. LNCS \textbf{6197}, Springer, Heidelberg, 2010, pp. 96--109.

\bibitem{tors}
A. Dujella, \emph{High rank elliptic curves with prescribed torsion}, \url{http://web.math.hr/~duje/tors/tors.html}

\bibitem{tors2}
A. Dujella, \emph{Infinite families of elliptic curves with high rank and prescribed torsion} \url{http://web.math.hr/~duje/tors/generic.html}

\bibitem{dm}
A. Dujella and M. Juki\'c Bokun, \emph{On the rank of elliptic curves over $\mathbf{Q}(i)$ with torsion group $\mathbf{Z}/4\mathbf{Z} \times \mathbf{Z}/4\mathbf{Z}$ }, Proc. Japan Acad. Ser. A Math. Sci. \textbf{86} (2010), 93--96.

\bibitem{jk4}
D.\ Jeon, C.~H.\ Kim, Y.\ Lee,
\emph{Families of elliptic curves over quartic number fields with prescribed torsion subgroups}, Math. Comp., to appear.

\bibitem{mir1}
M. Juki\'c Bokun, \emph{On the rank of elliptic curves over $\mathbf{Q}(\sqrt{-3})$ with torsion groups $\mathbf{Z}/3\mathbf{Z} \times \mathbf{Z}/3\mathbf{Z}$ and $\mathbf{Z}/3\mathbf{Z} \times \mathbf{Z}/6\mathbf{Z}$} , Proc. Japan Acad. Ser. A Math. Sci. \textbf{87} (2011), 61--64.

\bibitem{len}
H. W. Lenstra Jr., \emph{Factoring integers with elliptic curves}, Ann. of Math. \textbf{126} (1987) 649-
673.

\bibitem{Mont}
P.L. Montgomery, \emph{Speeding the Polard and elliptic curve methods of factorization},
Math. Comp. \textbf{48} (1987), 243--264.

\bibitem{Rab}
F. P. Rabarison, \emph{Torsion et rang des courbes elliptiques definies sur les
corps de nombres alg\'ebriques}, Doctorat de Universit\'e de Caen, 2008.

\bibitem{sil}
J.\ Silverman, \emph{The Arithmetic of Elliptic Curves}, Springer-Verlag, New York, 2009.


\end{thebibliography}
\end{document}